\documentclass[12pt,a4paper]{article}
\usepackage{amsmath, amssymb, theorem, latexsym, epsfig}
 \usepackage{graphics}

\newtheorem{theorem}{Theorem}[section]
\newtheorem{lemma}[theorem]{Lemma}
\newtheorem{proposition}[theorem]{Proposition}

\newtheorem{conjecture}[theorem]{Conjecture}
\newtheorem{corollary}[theorem]{Corollary}

\newtheorem{remark}[theorem]{Remark}

\allowdisplaybreaks

\def\neweq#1{\begin{equation}\label{#1}}
\def\endeq{\end{equation}}
\numberwithin{equation}{section}
\begin{document}
{\centering
\bfseries
{\LARGE A note on 5-cycle double covers}

\bigskip
\mdseries
{\large Arthur Hoffmann-Ostenhof}\footnote{supported by the FWF project P20543.}
\par
\upshape
\footnotesize \textit{Technical University of Vienna, Austria}\\

}

\noindent
\begin{abstract}
\noindent
The strong cycle double cover conjecture states that for every circuit $C$ of a bridgeless cubic graph 
$G$, there is a cycle double cover of $G$ which contains $C$. We conjecture that there is even a $5$-cycle double cover $S$ 
of $G$ which contains $C$, i.e. $C$ is a subgraph of one of the five $2$-regular subgraphs of $S$. We prove a necessary and sufficient condition 
for a $2$-regular subgraph to be contained in a $5$-cycle double cover of $G$. 
\end{abstract}


\vspace{0.5cm}

\noindent
Within the topic of cycle covers, a cycle is usually defined as a graph with even degree in every vertex. 
Since this is not the current standard graph theoretical definition of a cycle, we avoid the term cycle.  
We define a circuit to be a $2$-connected $2$-regular graph. For terminology not defined here we refer to \cite {Bo}. 
For a survey on cycle covers, see \cite {Z1}.

\noindent
A cycle double cover (CDC) of a graph $G$ with vertex degree at most $3$ is a set $\textbf{S}$ of $2$-regular subgraphs of $G$ 
such that every edge of $G$ is contained in exactly two elements of $\textbf{S}$.

\noindent
The cycle double cover conjecture (CDCC) states that every bridgeless cubic graph has a CDC.

\noindent
If a CDC $\textbf{S}$ of a graph $G$ satisfies $|\textbf{S}|=k$, then $\textbf{S}$ is called a $k$-CDC of $G$. 
There are several variations of the CDCC, as the strong-CDCC \cite {L} and the 5-CDCC \cite {C,P}.

\noindent
Strong-CDCC: For every given circuit $C$ of a bridgeless cubic graph $G$, there is a CDC $\textbf{S}$ with $C \in \textbf{S}$. 

\noindent
$5$-CDCC: Every bridgeless cubic graph has a $5$-CDC.

\noindent
For the subsequent proofs we need the following two lemmas which are stated in \cite {K,Z1} for cubic graphs but which 
hold as well for subdivisions of cubic graphs.

\begin{lemma}\label{tait2}
A subdivision of a cubic graph has a $4$-CDC if and only if it has a nowhere-zero $4$-flow.
\end{lemma}

\begin{lemma}\label{tait1}
Let $G$ be a subdivision of a cubic graph admitting a nowhere zero $4$-flow and $C'$ be a $2$-regular subgraph of $G$. Then there is a $4$-CDC $\textbf{S}$ of $G$ with $C' \in \textbf{S}$. 
\end{lemma}

\noindent
Let $C_0$ be a $2$-regular subgraph of a cubic graph $G$. 
We say $C_0$ is \textit{contained} in a CDC $\textbf{S}$ of $G$ if and only if $C_0$ is a subgraph of one of the elements of $\textbf{S}$.

\begin{proposition}\label{th1}
Let $G$ be a cubic graph with a $2$-regular subgraph $C_0 \subseteq G$. Moreover, let $\textbf{C} :=\{C_1,C_2,...,C_k\}$ be a set of 
$k$ $2$-regular subgraphs of $G$ with $C_0 \subseteq C_1$ such that

\noindent
1. every edge of $G$ is contained in at most two elements of $\textbf{C}$;

\noindent
2. the edges of $G$ which are contained in two elements of $\textbf{C}$ form a matching $M$ of $G$;

\noindent
3. $G-M$ has a nowhere-zero $4$-flow;

\noindent
then $G$ has a $(k+3)$-CDC which contains $C_0$.

\end{proposition}

\noindent
Proof: Suppose $G$ has such a set \textbf{C} of $2$-regular subgraphs as described above. Set $G':=G-M$ and let $C'$ be the $2$-regular subgraph which is constructed 
from $\textbf{C}$ by taking the symmetric difference of all elements of $\textbf{C}$. By condition 2 above, $C'$ is a $2$-regular subgraph of $G'$. Since $G'$ has a 
nowhere zero $4$-flow we can apply Lemma \ref{tait1}. Let $\textbf{S'}:=\{C',C_1',C_2',C_3'\}$ be a $4$-CDC of $G'$. 
Set $\textbf{S}:=\{C_1', C'_2, C'_3, C_1, C_2,...,C_k\}$. We claim that $\textbf{S}$ is a CDC of $G$. If $e \in M \subseteq E(G)$, then $e$ is 
covered twice by the elements of $\textbf{S} \cap \textbf{C}$. For $e \in E(G) - M$ it follows that $e \in C_m'$, $m \in \{1,2,3\}$. 
Consequently, there are two cases to consider:

\noindent
Case 1: $e \in E(C_i), i \in \{1,2,...,k\}$.\\
Then $e$ is covered precisely by $C_i$ and $C_m'$.

\noindent
Case 2: $e \notin E(C_i), i = 1,2,...,k$.\\
Then $e \notin C'$ and is thus covered precisely by $C_m'$ and $C_r'$ with $r \not= m$ and $r \in \{1,2,3\}$.

\noindent
Hence $\textbf{S}$ is a $(k+3)$-CDC of $G$ which contains $C_0$ since $C_0 \subseteq C_1$ by assumption. 

\newpage

\begin{theorem}\label{th2}
Let $G$ be a cubic graph with a $2$-regular subgraph $C_0 \subseteq G$. Then $G$ has a $5$-CDC which contains $C_0$ if and only if
 
\noindent
1. $G$ contains a matching $M$ such that $G-M$ has a nowhere-zero $4$-flow, and 

\noindent 
2. $G$ contains two $2$-regular subgraphs $C_1$, $C_2$ with $M=E(C_1) \cap E(C_2)$ and $C_0 \subseteq C_1$. 
\end{theorem}

\noindent
Proof: Suppose $M$, $C_1$, $C_2$ exist with the properties described in 1. and 2. above, then set $\textbf{C}:= \{C_1,C_2\}$ and 
apply Proposition \ref {th1}. Hence $G$ has a 5-CDC containing $C_0$.

\noindent
Conversely, suppose $G$ has a $5$-CDC $\textbf{S}:=\{C_1,C_2,...,C_5\}$ with $C_0 \subseteq C_1$. Suppose without loss of generality 
that $C_1 \cap C_2 \not = \emptyset $. Then $C_1 \cap C_2$ is a matching $M$ in $G$ since $\textbf{S}$ is a CDC and thus 
condition 2 is satisfied.\\ Moreover $G-M$ has a 4-CDC $\textbf{S'}:=\{C_1 \triangle C_2,C_3,C_4, C_5\}$ and thus $G-M$ has a 
nowhere-zero $4$-flow by Lemma \ref{tait2} which finishes the proof. 

\vspace{0.1cm}

\begin{remark}
Theorem \ref {th2} is useful for showing that for a given circuit $C$ in a snark $G$ there is a CDC of $G$ which contains $C$.
It follows by this theorem that it is sufficient to find a circuit $C'$ which intersects $C$ in a matching $M$ such that $G-M$ has a nowhere-zero $4$-flow.
If $G$ is the Petersen Graph, it suffices to show that $G-M$ is bridgeless. Thus, it is straightforward to see that every circuit of 
the Petersen Graph is contained in a $5$-CDC. 
\end{remark}

\begin{corollary}
A cubic graph $G$ has a $5$-CDC if and only if $G$ has a matching $M$ such that 
$G-M$ has a nowhere-zero $4$-flow and $G$ contains two $2$-regular subgraphs $C_1$ and $C_2$ with $E(C_1) \cap E(C_2)=M$.
\end{corollary}

\noindent
The preceding considerations lead us to formulate the following conjecture which we call the strong 5-CDCC.
This conjecture is a combination and strengthening of the 5-CDCC and the strong-CDCC.

\begin{conjecture}
For every given circuit $C$ of a bridgeless cubic graph $G$, there is a 5-CDC $\textbf{S}$ which contains $C$. 
\end{conjecture}

\noindent
Note that the strong 5-CDCC has been verified for all snarks of order less than $36$; see \cite {Z2}. 

\footnotesize
\bibliographystyle{plain}

\end{document}